# Random systems of polynomial equations. The expected number of roots under smooth analysis

DIEGO ARMENTANO[*] and MARIO WSCHEBOR[**]

*Centro de Matemática, Facultad de Ciencias, Universidad de la República, Calle Igua 4225, 11400 Montevideo, Uruguay. E-mail:* [*]*diego@cmat.edu.uy;* [**]*wschebor@cmat.edu.uy*

We consider random systems of equations over the reals, with $m$ equations and $m$ unknowns $P_i(t) + X_i(t) = 0$, $t \in \mathbb{R}^m$, $i = 1, \ldots, m$, where the $P_i$'s are non-random polynomials having degrees $d_i$'s (the "signal") and the $X_i$'s (the "noise") are independent real-valued Gaussian centered random polynomial fields defined on $\mathbb{R}^m$, with a probability law satisfying some invariance properties.

For each $i$, $P_i$ and $X_i$ have degree $d_i$.

The problem is the behavior of the number of roots for large $m$. We prove that under specified conditions on the relation signal over noise, which imply that in a certain sense this relation is neither too large nor too small, it follows that the quotient between the expected value of the number of roots of the perturbed system and the expected value corresponding to the centered system (i.e., $P_i$ identically zero for all $i = 1, \ldots, m$), tends to zero geometrically fast as $m$ tends to infinity. In particular, this means that the behavior of this expected value is governed by the noise part.

## 1. Introduction and main result

Let $f = (f_1, \ldots, f_m)$,

$$f_i(t) := \sum_{\|j\| \le d_i} a_j^{(i)} t^j \qquad (i = 1, \ldots, m) \tag{1}$$

be a system of $m$ polynomials in $m$ real variables. The notation in (1) is the following: $t := (t_1, \ldots, t_m)$ denotes a point in $\mathbb{R}^m$, $j := (j_1, \ldots, j_m)$ a multi-index of non-negative integers, $\|j\| = \sum_{h=1}^m j_h$, $t^j = t^{j_1} \cdots t^{j_m}$, $a_j^{(i)} = a_{j_1, \ldots, j_m}^{(i)}$. $d_i$ is the degree of the polynomial $f_i$.

We denote by $N^f(V)$ the number of roots of the system of equations

$$f_i(t) = 0 \qquad (i=1,\ldots,m)$$

lying in the subset $V$ of $\mathbb{R}^m$. Also $N^f = N^f(\mathbb{R}^m)$.

If we choose at random the coefficients $\{a_j^{(i)}\}$, $N^f(V)$ becomes a random variable. Classical results in the case of one polynomial in one variable seem to have started with the work of Marc Kac [6] (see the book by Bharucha-Reid and Sambandham [3]). Here we will be interested in systems with $m > 1$, and more specifically, in large values of $m$. This appears to be of a quite different nature than the case $m = 1$ and generally speaking, little is known on the distribution of the random variable $N^f(V)$ (or $N^f$) even for simple choices of the probability law on the coefficients.

In 1992, Shub and Smale [10] computed the expectation of $N^f$ when the coefficients are Gaussian centered independent random variables having variances:

$$E[(a_j^{(i)})^2] = \frac{d_i!}{j_1! \cdots j_m! (d_i - \|j\|)!}. \qquad (2)$$

Their result was

$$E(N^f) = (d_1 \cdots d_m)^{1/2}. \qquad (3)$$

Some extensions of their work, including new results for one polynomial in one variable, can be found in [5]. There are also other extensions to multi-homogeneous systems in [8], and, partially, to sparse systems in [7] and [9]. A similar question for the number of critical points of real-valued polynomial random functions has been considered in a recent paper by Dedieu and Malajovich [4].

A general formula for $E(N^f(V))$ when the random functions $f_i$ $(i=1,\ldots,m)$ are stochastically independent and their law is centered and invariant under the orthogonal group on $\mathbb{R}^m$ can be found in [1]. This includes the Shub–Smale formula (3) as a special case. Very little in known on higher moments. The only published results of which the authors are aware concern asymptotic variances as $m \to +\infty$ (see [1] for non-polynomial systems and [11] for the Kostlan–Shub–Smale model).

The aim of this paper is to remove the hypothesis that the coefficients have zero expectation (in some cases, this has been considered for one polynomial in one variable in the above-mentioned paper [5]).

One way to look at this problem is to start with a non-random system

$$P_i(t) = 0 \qquad (i=1,\ldots,m), \qquad (4)$$

perturb it with a polynomial noise $\{X_i(t) : i = 1,\ldots,m\}$, that is, consider

$$P_i(t) + X_i(t) = 0 \qquad (i=1,\ldots,m)$$

and ask what one can say about the number of roots of the new system. Of course, to obtain results on $E(N^{P+X})$ we need a certain number of hypotheses both on the "noise" $X$ and the class of polynomial "signals" $P$, especially the relation between the size of $P$ and the probability distribution of $X$.

Roughly speaking, we prove in Theorem 2 that if the relation signal over noise is neither too big nor too small, in a sense that will be made precise later on, there exist positive constants $C, \theta$, $0 < \theta < 1$ such that

$$E(N^{P+X}) \leq C\theta^m E(N^X). \tag{5}$$

Inequality (5) becomes of interest if the starting non-random system (4) has a large number of roots, possibly infinite, and $m$ is large. In this situation, the effect of adding polynomial noise is a reduction at a geometric rate of the expected number of roots, as compared to the centered case in which all the $P_i$'s are identically zero.

Notice that in formula (5), $E(N^X)$ can be computed by means of a nice formula once we know the probability distribution of the noise and that the constants $C, \theta$ can be explicitly estimated from the hypotheses.

We will assume throughout that the polynomial noise $X$ is Gaussian and centered, the real-valued random processes

$$X_1(\cdot), \ldots, X_m(\cdot)$$

defined on $\mathbb{R}^m$ are independent, with covariance functions

$$R^{X_i}(s,t) := E(X_i(s)X_i(t)) \qquad (i = 1, \ldots, m)$$

depending only on the scalar product $\langle s, t \rangle$, that is: $R^{X_i}(s,t) = Q^{(i)}(\langle s, t \rangle)$, where

$$Q^{(i)}(u) = \sum_{k=0}^{d_i} c_k^{(i)} u^k, \qquad u \in \mathbb{R} \ (i = 1, \ldots, m). \tag{6}$$

In this case, it is known that a necessary and sufficient condition for $Q^{(i)}(\langle s, t \rangle)$ to be a covariance is that $c_k^{(i)} \geq 0$ for all $k = 0, \ldots, d_i$ and the process $X_i$ can be written as

$$X_i(t) = \sum_{\|j\| \leq d_i} a_j^{(i)} t^j,$$

where the random variables $a_j^{(i)}$ are centered Gaussian, independent and

$$\mathrm{Var}(a_j^{(i)}) = c_{\|j\|}^{(i)} \frac{\|j\|!}{j!} \qquad (i = 1, \ldots, m; \|j\| \leq d_i)$$

(for a proof, see, e.g., [1]).

The Shub–Smale model (2) corresponds to the particular choice

$$c_k^{(i)} = \binom{d_i}{k} \qquad (k = 0, 1, \ldots, d_i)$$

which implies

$$Q^{(i)}(u) = (1+u)^{d_i} \qquad (i=1,\ldots,m). \tag{7}$$

We will use the following notations:

$Q_u^{(i)}, Q_{uu}^{(i)}$ denote the successive derivatives of $Q^{(i)}$. We assume that $Q^{(i)}(u), Q_u^{(i)}(u)$ do not vanish for $u \geq 0$. Put, for $x \geq 0$:

$$q_i(x) := \frac{Q_u^{(i)}}{Q^{(i)}}, \tag{8a}$$

$$r_i(x) := \frac{Q^{(i)} Q_{uu}^{(i)} - (Q_u^{(i)})^2}{(Q^{(i)})^2}, \tag{8b}$$

$$h_i(x) := 1 + x \frac{r_i(x)}{q_i(x)}. \tag{8c}$$

In (8a) and (8b), the functions in the right-hand side are computed at the point $x$.

In [1] the following statement was proved:

**Theorem 1.** *For any Borel set $V \subset \mathbb{R}^m$ we have:*

$$E(N^X(V)) = \frac{1}{\sqrt{2}\pi^{(m+1)/2}} \Gamma\left(\frac{m}{2}\right) \int_V \left[\prod_{i=1}^m q_i(\|t\|^2)\right]^{1/2} \cdot E_h(\|t\|^2) \, dt \tag{9}$$

*where*

$$E_h(x) = E\left(\left[\sum_{i=1}^m h_i(x) \xi_i^2\right]^{1/2}\right)$$

*and $\xi_1, \ldots, \xi_m$ denote independent standard normal random variables.*

*Remark.* In fact, Theorem 1 is a special case of a general theorem (see [1]), in which the covariance function of the random field is invariant under the action of the orthogonal group, and not only a function of the scalar product.

Before the statement of our main result, Theorem 2 below, we need to introduce some additional notations and hypotheses.

We will assume that each polynomial $Q^{(i)}$ does not vanish for $u \geq 0$, which amounts to saying that for each $t$ the (one-dimensional) distribution of $X_i(t)$ does not degenerate. Also, $Q^{(i)}$ has effective degree $d_i$, that is,

$$c_{d_i}^{(i)} > 0 \qquad (i=1,\ldots,m).$$

An elementary calculation then shows that for each polynomial $Q^{(i)}$, as $u \to +\infty$:

$$q_i(u) \sim \frac{d_i}{1+u}, \tag{10a}$$

$$h_i(u) \sim \frac{c^{(i)}_{d_i-1}}{d_i c^{(i)}_{d_i}} \cdot \frac{1}{1+u}. \tag{10b}$$

Since we are interested in the large $m$ asymptotic and the polynomials $P, Q$ can vary with $m$, we will require somewhat more than relations (10a) and (10b), as specified in the following hypotheses:

(H$_1$) $h_i$ is independent of $i$ ($i = 1, \ldots, m$) (but may vary with $m$). We put $h = h_i$. Of course, if the polynomials $Q^{(i)}$ do not depend on $i$, this hypothesis is satisfied. But there are more general cases, such as covariances having the form $Q(u)^{l_i}$ ($i = 1, \ldots, m$).

(H$_2$) There exist positive constants $D_i, E_i$ ($i = 1, \ldots, m$) and $\underline{q}$ such that

$$0 \leq D_i - (1+u)q_i(u) \leq \frac{E_i}{1+u} \quad \text{and} \quad (1+u)q_i(u) \geq \underline{q} \tag{11}$$

for all $u \geq 0$, and moreover

$$\max_{1 \leq i \leq m} D_i, \qquad \max_{1 \leq i \leq m} E_i$$

are bounded by constants $\overline{D}, \overline{E}$, respectively, which are independent of $m$; $\underline{q}$ is also independent of $m$.

Also, there exist positive constants $\underline{h}, \overline{h}$ such that

$$\underline{h} \leq (1+u)h(u) \leq \overline{h} \tag{12}$$

for $u \geq 0$.

Notice that the auxiliary functions $q_i, r_i, h$ ($i = 1, \ldots, m$) will also vary with $m$. To simplify somewhat the notation, we are dropping the parameter $m$ in $P, Q, q_i, r_i, h$. However, in (H$_2$) the constants $\underline{h}, \overline{h}$ do not depend on $m$ (see the examples after the statement of Theorem 2 below).

With respect to (H$_2$), it is clear that for each $i$, $q_i$ will satisfy (11) with the possible exception of the first inequality, and $(1+u)h(u) \leq \overline{h}$ for some positive constant $\overline{h}$, from the definitions (8a), (8c), (10a), (10b) and the conditions on the coefficients of $Q^{(i)}$. However, it is not self-evident from the definition (8c) that $h(u) \geq 0$ for $u \geq 0$. This will become clear in the proof of Theorem 2 below.

A second set of hypotheses on the system concerns the relation between the "signal" $P$ and the "noise" $X$, which roughly speaking should neither be too small nor too big.

Let $P$ be a polynomial in $m$ real variables with real coefficients having degree $d$ and $Q$ a polynomial in one variable with non-negative coefficients, also having degree $d$,

$Q(u) = \sum_{k=0}^{d} c_k u^k$. We assume that $Q$ does not vanish on $u \geq 0$ and $c_d > 0$. Define

$$H(P,Q) := \sup_{t \in \mathbb{R}^m} \left\{ (1+\|t\|) \cdot \left\| \nabla\left(\frac{P}{\sqrt{Q(\|t\|^2)}}\right)(t) \right\| \right\},$$

$$K(P,Q) := \sup_{t \in \mathbb{R}^m \setminus \{0\}} \left\{ (1+\|t\|^2) \cdot \left| \frac{\partial}{\partial \rho}\left(\frac{P}{\sqrt{Q(\|t\|^2)}}\right)(t) \right| \right\},$$

where $\frac{\partial}{\partial \rho}$ denotes the derivative in the direction defined by $\frac{t}{\|t\|}$, at each point $t \neq 0$.

For $r > 0$, put:

$$L(P,Q,r) := \inf_{\|t\| \geq r} \frac{P(t)^2}{Q(\|t\|^2)}.$$

One can check by means of elementary computations that for each pair $P,Q$ as above, one has

$$H(P,Q) < \infty, \qquad K(P,Q) < \infty.$$

With these notations, we introduce the following hypotheses on the systems $P,Q$, as $m$ grows:

(H$_3$)

$$A_m = \frac{1}{m} \cdot \sum_{i=1}^{m} \frac{H^2(P_i, Q^{(i)})}{i} = o(1) \qquad \text{as } m \to +\infty, \tag{13a}$$

$$B_m = \frac{1}{m} \cdot \sum_{i=1}^{m} \frac{K^2(P_i, Q^{(i)})}{i} = o(1) \qquad \text{as } m \to +\infty. \tag{13b}$$

(H$_4$) There exist positive constants $r_0, \ell$ such that if $r \geq r_0$:

$$L(P_i, Q^{(i)}, r) \geq \ell \qquad \text{for all } i = 1, \ldots, m.$$

**Theorem 2.** *Under the hypotheses* (H$_1$), $\ldots,$ (H$_4$), *one has*

$$E(N^{P+X}) \leq C\theta^m E(N^X), \tag{14}$$

*where $C, \theta$ are positive constants, $0 < \theta < 1$.*

## 1.1. Remarks on the statement of Theorem 2

1. In fact, we will see in the proof of the theorem how one can get explicitly from the hypotheses first the value of $\theta$ and then $m_0$ and the constant $C$ in such a way that whenever $m \geq m_0$, inequality (14) holds true.

    A possible choice is as follows:

- Choose $r_0$ from (H$_4$),

$$\theta_1 = \max\left\{\frac{r_0}{\sqrt{r_0^2+1/2}}, e^{-\ell/2}\right\}, \qquad \theta = \frac{1+\theta_1}{2}.$$

- Let us put $F_i = E_i/D_i$ $(i=1,\ldots,m)$ and $\bar{F} = \max\{F_1,\ldots,F_m\}$. From the hypotheses, one has $\bar{F} \leq \bar{E}/\underline{q}$. Let $\tau > 0$ such that:

$$\frac{\bar{F}}{1+\tau^2 r_0^2} < \frac{1}{2}\frac{1}{1+r_0^2}. \tag{15}$$

Choose $m_0$ (using (H$_3$)) so that if $m \geq m_0$ one has:

$$\begin{aligned} e^{[mA_m/\underline{q}+mB_m/(\underline{h}\underline{q})]/2}\theta_1^m \sqrt{m} &\leq \theta^m, \\ \pi\left(\frac{\tau^2 r_0^2}{1+\tau^2 r_0^2}\right)^{(m-1)/2} &< \frac{e^{-2}}{\sqrt{m}}. \end{aligned} \tag{16}$$

*Then*, (14) is satisfied for $m \geq m_0$, with

$$C = 30 \cdot \frac{\bar{h}}{\underline{h}} \frac{\sqrt{1+r_0^2}}{r_0}. \tag{17}$$

2. It is obvious that our problem does not depend on the order in which the equations

$$P_i(t) + X_i(t) = 0 \quad (i=1,\ldots,m)$$

appear. However, conditions (13a) and (13b) in hypothesis (H$_3$) do depend on the order. One can state them by saying that there exists an order $i=1,\ldots,m$ on the equations, such that (13a) and (13b) hold true.

3. Condition (H$_3$) can be interpreted as a bound on the quotient signal over noise. In fact, it concerns the gradient of this quotient. In (13b) the radial derivative appears, which happens to decrease faster as $\|t\| \to \infty$ than the other components of the gradient.

   Clearly, if $H(P_i,Q^{(i)}), K(P_i,Q^{(i)})$ are bounded by fixed constants, (13a) and (13b) are verified. Also, some of them may grow as $m \to +\infty$ provided (13a) and (13b) remain satisfied.

4. Hypothesis (H$_4$) goes – in some sense – in the opposite direction: For large values of $\|t\|$ we need a lower bound of the relation signal over noise.

5. A result of the type of Theorem 2 can not be obtained without putting some restrictions on the relation signal over noise. In fact, consider the system

$$P_i(t) + \sigma X_i(t) = 0 \qquad (i=1,\ldots,m), \tag{18}$$

where $\sigma$ is a positive real parameter. For generic $P$, as $\sigma \downarrow 0$ the expected value of the number of roots of (18) tends to the number of roots of $P_i(t) = 0$ $(i=1,\ldots,m)$.

In this case, the relation signal over noise tends to infinity. On the other hand, if we let $\sigma \to +\infty$, the relation signal over noise tends to zero and the expected number of roots will tend to $E(N^X)$.

## 2. Some examples

### 2.1. Shub–Smale

In the Shub–Smale model for the noise, $Q^{(i)}$ is given by (7). Then,

$$q_i(u) = \frac{d_i}{1+u}, \qquad h_i(u) = h(u) = \frac{1}{1+u}.$$

We assume that the degrees $d_i$ are uniformly bounded. So, (H$_1$) and (H$_2$) hold true. Of course, conditions (H$_3$) and (H$_4$) also depend on the signal.

We are going to give two simple examples. Let

$$P_i(t) = \|t\|^{d_i} - r^{d_i},$$

where $d_i$ is even and $r$ is positive and remains bounded as $m$ varies. One has:

$$\frac{\partial}{\partial \rho}\left(\frac{P_i}{\sqrt{Q^{(i)}}}\right)(t) = \frac{d_i \|t\|^{d_i-1} + d_i r^{d_i} \|t\|}{(1+\|t\|^2)^{d_i/2+1}} \leq \frac{d_i(1+r^{d_i})}{(1+\|t\|^2)^{3/2}},$$

$$\nabla\left(\frac{P_i}{\sqrt{Q^{(i)}}}\right)(t) = \frac{d_i \|t\|^{d_i-2} + d_i r^{d_i}}{(1+\|t\|^2)^{d_i/2+1}} \cdot t$$

which implies

$$\left\|\nabla\left(\frac{P_i}{\sqrt{Q^{(i)}}}\right)(t)\right\| \leq \frac{d_i(1+r^{d_i})}{(1+\|t\|^2)^{3/2}}.$$

Again, since the degrees $d_1, \ldots, d_m$ are bounded by a constant that does not depend on $m$, (H$_3$) follows. (H$_4$) also holds under the same hypothesis.

To illustrate a numerical example of the values of $\theta$, $C$ and $m_0$, let us assume that the radius $r = 1$. Then:

- We can choose $r_0 = 2$. It turns out that $\ell = 9/25$ if $\bar{D} \leq 4$ and $\ell = (2^{\bar{D}} - 1)^2/5^{\bar{D}}$ if $\bar{D} \geq 5$.
- $\theta = (3 + 2\sqrt{2})/6$.
- $C = 15\sqrt{5}$ and $m_0$ can be chosen so that if $m \geq m_0$, then

$$e^{c_1} m^{c_1+1/2} \leq \kappa^m, \qquad \text{where } c_1 = 8\bar{D}^2, \kappa = \frac{3}{4\sqrt{2}} + \frac{1}{2} > 1.$$

Notice that an interest in this choice of the $P_i$'s lies in the fact that obviously the system $P_i(t) = 0$ $(i = 1, \ldots, m)$ has infinite roots (all points in the sphere of radius $r$ centered at the origin are solutions), but the expected number of roots of the perturbed system is geometrically smaller than the Shub–Smale expectation, when $m$ is large.

Our second example is the following: Let $T$ be a polynomial of degree $d$ in one variable that has $d$ distinct real roots. Define:

$$P_i(t_1, \ldots, t_m) = T(t_i) \qquad (i = 1, \ldots, m).$$

One can easily check that the system verifies our hypotheses, so that there exist $C, \theta$ positive constants, $0 < \theta < 1$ such that

$$E(N^{P+X}) \leq C\theta^m d^{m/2},$$

where we have used the Shub–Smale formula when the degrees are all the same. On the other hand, it is clear that $N^P = d^m$ so that the diminishing effect of the noise on the number of roots can be observed. A number of variations of these examples for $P$ can be constructed, but we will not pursue the subject here.

## 2.2. $Q^{(i)} = Q$, only real roots

Assume all the $Q^{(i)}$ are equal, $Q^{(i)} = Q$ and $Q$ has only real roots. Since $Q$ does not vanish on $u \geq 0$, all the roots should be strictly negative, say

$$-\alpha_1, \ldots, -\alpha_d,$$

where

$$0 < \alpha_1 \leq \alpha_2 \leq \cdots \leq \alpha_d.$$

With no loss of generality, we may assume that $\alpha_1 \geq 1$. If this were not the case, we perform a homothecy of the space $\mathbb{R}^m$, centered at the origin and with factor equal to $\sqrt{\alpha_1}$. The number of roots remains unchanged, and the new $Q$ has $\alpha_1 \geq 1$.

We will assume again that the degree $d$ of $Q$ is bounded by a fixed constant $\overline{d}$ (one should take into account that $Q$ may vary with $m$), as well as the roots

$$\alpha_k \leq \overline{\alpha} \qquad (k = 1, \ldots, d)$$

for some constant $\overline{\alpha}$. A direct computation gives:

$$q_i(u) = q(u) = \sum_{k=1}^{d} \frac{1}{u + \alpha_k}, \qquad h_i(u) = h(u) = \frac{1}{q(u)} \sum_{k=1}^{d} \frac{\alpha_k}{(u + \alpha_k)^2}.$$

One verifies (11), choosing $D_i = d$, $E_i = d(\alpha_d - 1)$. Similarly, a direct computation gives (12).

Again let us consider the particular example of signals:

$$P_i(t) = \|t\|^{d_i} - r^{d_i},$$

where $d_i$ is even and for each $i = 1, \ldots, m$ and $r$ is positive and remains bounded as $m$ varies.

$$\left| \frac{\partial}{\partial \rho} \left( \frac{P_i}{\sqrt{Q^{(i)}}} \right) \right| \leq d_i(\overline{\alpha} + r^{d_i}) \frac{1}{(1 + \|t\|^2)^{3/2}}$$

so that $K(P_i, Q^{(i)})$ is uniformly bounded. A similar computation shows that $H(P_i, Q^{(i)})$ is uniformly bounded. Finally, it is obvious that

$$L(P_i, Q^{(i)}, r) \geq \left( \frac{1}{1 + \overline{\alpha}} \right)^{\overline{d}}$$

for $i = 1, \ldots, m$ and any $r \geq 1$. So the conclusion of Theorem 2 can be applied.

One can similarly check that the second polynomial system in the previous example also works with respect to this noise.

## 2.3. More general examples

Assume that the noise has covariance with the form

$$Q^{(i)}(u, v, w) = [Q(u)]^{l_i} \qquad (i = 1, \ldots, m),$$

where $Q$ is a polynomial in one variable having degree $\nu$ with positive coefficients, $Q(u) = \sum_{k=0}^{\nu} b_k u^k$. $Q$ may depend on $m$, as well as the exponents $l_1, \ldots, l_m$. Notice that $d_i = \nu \cdot l_i$ ($i = 1, \ldots, m$).

The Shub–Smale case corresponds to the simple choice $Q(u) = 1 + u$, $l_i = d_i$ ($i = 1, \ldots, m$).

One has:

$$q_i(u) = l_i \frac{Q'(u)}{Q(u)},$$

$$h_i(u) = h(u) = 1 - u \frac{Q'^2(u) - Q(u)Q''(u)}{Q(u)Q'(u)}$$

so that (H$_1$) is satisfied.

We will require the coefficients $b_0, \ldots, b_\nu$ of the polynomial $Q$ to verify the conditions

$$b_k \leq \frac{\nu - k + 1}{k} b_{k-1} \qquad (k = 1, 2, \ldots, \nu).$$

Moreover, we assume that

$$l_1, \ldots, l_m, \nu$$

are bounded by a constant independent of $m$ and there exist positive constant $\underline{b},\overline{b}$ such that

$$\underline{b} \leq b_0, b_1, \ldots, b_\nu \leq \overline{b}.$$

Under these conditions, one can check that ($H_2$) holds true, with $D_i = d_i$ $(i = 1, \ldots, m)$.

For the relation signal over noise, conditions are similar to the previous example.

Notice that already if $\nu = 2$ and we choose for $Q$ the fixed polynomial:

$$Q(u) = 1 + 2au + bu^2$$

with $0 < a \leq 1$, $\sqrt{b} > a \geq b > 0$, then the conditions in this example are satisfied, but the polynomial $Q$ (hence $Q^{d_i}$) does not have real roots, so that it is not included in Example 2.2.

## 3. Proof of Theorem 2

**Proof of Theorem 2.** Let

$$Z_j(t) = \frac{P_j(t) + X_j(t)}{\sqrt{Q^{(j)}(\|t\|^2)}} \quad (j = 1, \ldots, m)$$

and

$$Z = \begin{pmatrix} Z_1 \\ \vdots \\ Z_m \end{pmatrix}.$$

Clearly,

$$N^{P+X}(V) = N^Z(V)$$

for any subset $V$ of $\mathbb{R}^m$.

Clearly, the Gaussian random fields $\{Z_j(t) : t \in \mathbb{R}^m\}$ $(j = 1, \ldots, m)$ are independent and

$$\operatorname{Var}(Z_j^2(t)) = E(\tilde{Z}_j^2(t)) = 1 \tag{19}$$

for all $j = 1, \ldots, m$ and all $t \in \mathbb{R}^m$, where $\tilde{Z}_j(t) = Z_j(t) - E(Z_j(t))$.

Differentiating in (19) with respect to $t_\alpha (\alpha = 1, \ldots, m)$ we obtain that $E(\frac{\partial}{\partial t_\alpha}\tilde{Z}_j(t)\tilde{Z}_j(t)) = 0$, $(j = 1, \ldots, m)$. Since the joint distribution is Gaussian, this implies that $\nabla \tilde{Z}_j(t)$ is independent of $\tilde{Z}_j(t)$, that is, $Z'(t)$ and $Z(t)$ are independent.

We apply the Rice formula (see [2] for a complete proof) to compute $E(N^Z(V))$, that is:

$$E(N^Z(V)) = \int_V E(|\det(Z'(t))| \, | \, Z(t) = 0) \cdot p_{Z(t)}(0) \, dt,$$

where $p_\xi(\cdot)$ denotes the density of the probability distribution of the random vector $\xi$, whenever it exists and $f'(t):\mathbb{R}^m \to \mathbb{R}^m$ the derivative of the function $f:\mathbb{R}^m \to \mathbb{R}^m$ at the point $t$. Because of the independence between $Z'(t)$ and $Z(t)$, we can erase the condition in the conditional expectation, obtaining:

$$E(N^Z(V)) = \int_V E(|\det(Z'(t))|) \cdot \frac{1}{(2\pi)^{m/2}} \\ \times e^{[-(P_1(t)^2/Q^{(1)}(\|t\|^2)+\cdots+P_m(t)^2/Q^{(m)}(\|t\|^2))/2]}\,dt. \quad (20)$$

Our next problem is the evaluation of $E(|\det(Z'(t))|)$.

A direct computation of covariances gives:

$$\mathrm{Cov}\left(\frac{\partial Z_i}{\partial t_\alpha}(t), \frac{\partial Z_j}{\partial t_\beta}(t)\right) = \delta_{ij}\frac{\partial^2}{\partial s_\alpha \partial t_\beta}R^{\tilde{Z}_i}(s,t)|_{s=t}$$
$$= \delta_{ij}[r_i(\|t\|^2)t_\alpha t_\beta + q_i(\|t\|^2)\delta_{\alpha\beta}]$$

for $i,j,\alpha,\beta = 1,\ldots,m$, where the functions $q_i, r_i$ have been defined in (8a) and (8b).

For each $t \neq 0$, let $U_t$ be an orthogonal transformation of $\mathbb{R}^m$ that takes the first element of the canonical basis into the unit vector $\frac{t}{\|t\|}$. Then

$$\mathrm{Var}\left(\frac{U_t \nabla Z_j(t)}{\sqrt{q_j(\|t\|^2)}}\right) = \mathrm{Diag}(h(\|t\|^2), 1, \ldots, 1), \quad (21)$$

where we denote the gradient $\nabla Z_j(t)$ as a column vector,

$$\nabla Z_j(t) = \begin{pmatrix} \frac{\partial Z_j}{\partial t_1}(t) \\ \vdots \\ \frac{\partial Z_j}{\partial t_m}(t) \end{pmatrix}.$$

$\mathrm{Diag}(\lambda_1,\ldots,\lambda_m)$ denotes the $m\times m$ diagonal matrix with elements $\lambda_1,\ldots,\lambda_m$ in the diagonal and the function $h$ has been defined in (8c). So we can write

$$\frac{U_t \nabla Z_j(t)}{\sqrt{q_j(\|t\|^2)}} = \zeta_j + \alpha_j \qquad (j=1,\ldots,m),$$

where $\zeta_j$ is a Gaussian centered random vector in $\mathbb{R}^m$ having covariance given by (21), $\zeta_1,\ldots,\zeta_m$ are independent and $\alpha_j$ is the non-random vector

$$\alpha_j = \frac{U_t \nabla(P_j(t)/\sqrt{Q^{(j)}(\|t\|^2)})}{\sqrt{q_j(\|t\|^2)}} = \begin{pmatrix} \alpha_{1j} \\ \vdots \\ \alpha_{mj} \end{pmatrix} \qquad (j=1,\ldots,m). \quad (22)$$

Denote $T$ as the $m \times m$ random matrix having columns $\eta_j = \zeta_j + \alpha_j$ $(j = 1, \ldots, m)$. We have

$$|\det(Z'(t))| = |\det(T)| \cdot \prod_{i=1}^{m}(q_i(\|t\|^2))^{1/2}$$

so that

$$E(|\det(Z'(t))|) = E(|\det(T)|) \cdot \prod_{i=1}^{m}(q_i(\|t\|^2))^{1/2}. \tag{23}$$

The components $\zeta_{ij}$ $(i = 1, \ldots, m)$ of $\zeta_j$ are Gaussian centered independent and

$$\mathrm{Var}(\zeta_{ij}) = 1 \qquad \text{for } i = 2, \ldots, m; j = 1, \ldots, m$$
$$\mathrm{Var}(\zeta_{1j}) = h(\|t\|^2) \quad \text{for } j = 1, \ldots, m.$$

Put

$$\tilde{\alpha}_j = \begin{pmatrix} \alpha_{1j}/\sqrt{h(\|t\|^2)} \\ \alpha_{2j} \\ \vdots \\ \alpha_{mj} \end{pmatrix}$$

so that

$$|\det(T)| = \sqrt{h(\|t\|^2)} \cdot |\det(\tilde{T})|,$$

where $\tilde{T}$ is the random matrix having columns $\tilde{\eta}_j = \psi_j + \tilde{\alpha}_j$ with $\psi_1, \ldots, \psi_m$ i.i.d. standard normal in $\mathbb{R}^m$.

$|\det(\tilde{T})|$ is the volume of the parallelotope generated by $\tilde{\eta}_1, \ldots, \tilde{\eta}_m$, that is, the set of vectors in $\mathbb{R}^m$ which can be expressed as

$$\sum_{j=1}^{m} a_j \tilde{\eta}_j$$

with $0 \leq a_i \leq 1$ $(i = 1, \ldots, m)$. Hence,

$$|\det(T)| = \sqrt{h(\|t\|^2)} \cdot \|\tilde{\eta}_m\| \cdot \prod_{j=1}^{m-1} d(\tilde{\eta}_j, S_j), \tag{24}$$

where $S_j$ $(j = 1, \ldots, m-1)$ stands for the subspace of $\mathbb{R}^m$ generated by $\tilde{\eta}_{j+1}, \ldots, \tilde{\eta}_m$ and $d$ is Euclidean distance. Notice that for $j = 1, \ldots, m-1$, using invariance of the standard Gaussian distribution under the orthogonal group,

$$E(d(\tilde{\eta}_j, S_j) \mid \tilde{\eta}_{j+1}, \ldots, \tilde{\eta}_m) = E(\|\xi_j + \tilde{a}_j\|_j), \tag{25}$$

where $\|\cdot\|_j$ denotes Euclidean norm in $\mathbb{R}^j$ ($\|\cdot\| = \|\cdot\|_m$), $\xi_j$ is standard normal in $\mathbb{R}^j$ and $\tilde{a}_j$ is the orthogonal projection of $\tilde{\alpha}_j$ onto the orthogonal complement of $S_j$ in $\mathbb{R}^m$ (which is identified here with $\mathbb{R}^j$), so that $\|\tilde{a}_j\|_j \leq \|\tilde{\alpha}_j\|$. Introduce now the function

$$\gamma_j(c) = E(\|\xi_j + c\|_j),$$

where $c \in \mathbb{R}^j$ is non-random. It is clear that $\gamma_j$ is a function only of $\|c\|_j$, which is in fact increasing and

$$\gamma_j(c) \leq \left(1 + \|c\|_j^2 \frac{1}{2j}\right)\gamma_j(0).$$

(See the auxiliary Lemma 1 after this proof.) So, if we denote by $c_j$ a non-random vector in $\mathbb{R}^j$ such that $\|c_j\|_j = \|\tilde{\alpha}_j\|_j$, it follows from (24) and (25), by successive conditioning that

$$E(|\det T|) \leq \sqrt{h(\|t\|^2)} \cdot \prod_{j=1}^{m} E(\|\xi_j + c_j\|_j) \leq \sqrt{h(\|t\|^2)} \cdot \left(\prod_{j=1}^{m} \gamma_j(0)\right)\left(\prod_{j=1}^{m}\left(1 + \|c\|_j^2 \frac{1}{2j}\right)\right).$$

Using (20) and (23) we get:

$$E(N^Z) \leq \frac{1}{(2\pi)^{m/2}} L_m \cdot \int_{\mathbb{R}^m} \left\{\sqrt{h(\|t\|^2)} \cdot \left(\prod_{i=1}^{m} q_i(\|t\|^2)\right)^{1/2}\right.$$
$$\left.\times \exp\left[-\frac{1}{2}\sum_{i=1}^{m}\frac{P_i(t)^2}{Q^{(i)}(\|t\|^2)} + \frac{1}{2}\sum_{j=1}^{m}\|c_j\|_j^2 \frac{1}{j}\right]\right\} dt, \quad (26)$$

where

$$L_m = \prod_{j=1}^{m} E(\|\xi_j\|_j) = \frac{1}{\sqrt{2\pi}} 2^{(m+1)/2} \Gamma\left(\frac{m+1}{2}\right).$$

Our final task is to obtain an adequate bound for the integral in (26). For $j = 1, \ldots, m$ (use $(H_2)$):

$$|\tilde{\alpha}_{1j}| = \frac{1}{\sqrt{h(\|t\|^2)q_j(\|t\|^2)}} \cdot \left|\frac{\partial}{\partial \rho}\frac{P_j(\|t\|^2)}{\sqrt{Q^{(j)}(\|t\|^2)}}\right| \leq \frac{1}{\sqrt{\underline{hq}}} K(P_j, Q^{(j)})$$

and

$$\|\alpha_j\| = \frac{\|\nabla(P_j(t)/\sqrt{Q^{(j)}(\|t\|^2)})\|}{\sqrt{q_j(\|t\|^2)}} \leq \frac{1}{\sqrt{\underline{q}}} H(P_j, Q^{(j)}).$$

Then, if we bound $\|\tilde{\alpha}_j\|^2$ by:

$$\|\tilde{\alpha}_j\|^2 \leq |\tilde{\alpha}_{1j}|^2 + \|\alpha_j\|^2$$

we obtain
$$\|\tilde{\alpha}_j\|^2 \leq \frac{1}{\underline{h}\underline{q}} K^2(P_j, Q^{(j)}) + \frac{1}{\underline{q}} H^2(P_j, Q^{(j)}),$$

which implies (using (H$_3$))
$$\sum_{j=1}^{m} \|c_j\|_j^2 \cdot \frac{1}{j} \leq \frac{1}{\underline{q}} m A_m + \frac{1}{\underline{h}\underline{q}} m B_m.$$

Replacing in (26) we get the bound:
$$E(N^Z) \leq s_m H_m, \tag{27}$$

where
$$s_m = \left(\frac{\overline{h}}{\underline{h}}\right)^{1/2} \cdot e^{(mA_m/\underline{q} + mB_m/(\underline{h}\underline{q}))/2} \tag{28}$$

and
$$H_m = \frac{1}{\sqrt{2}\pi^{(m+1)/2}} \Gamma\left(\frac{m}{2}\right) \\ \times \int_{\mathbb{R}^m} \left(\prod_{i=1}^{m} q_i(\|t\|^2)\right)^{1/2} \sqrt{h(\|t\|^2)} E(\|\xi_m\|) e^{-(\sum_{i=1}^{m} P_i(t)^2/Q^{(i)}(\|t\|^2))/2} \, dt. \tag{29}$$

The integrand in (29) is the same as in formula (9) giving the expectation in the centered case, except for the exponential, which will help for large values of $\|t\|$.

Let us write $H_m$ as
$$H_m = H_m^{(1)}(r) + H_m^{(2)}(r),$$

where $H_m^{(1)}(r)$ corresponds to integrating on $\|t\| \leq r$ and $H_m^{(2)}(r)$ on $\|t\| > r$ instead of the whole $\mathbb{R}^m$ in formula (29). We first choose $r = r_0$ so that, using $H_4$:
$$H_m^{(2)}(r_0) \leq e^{-\ell \, m/2} E(N^X). \tag{30}$$

We now turn to $H_m^{(1)}(r)$. We have, bounding the exponential in the integrand by 1 and using hypothesis (H$_2$):
$$H_m^{(1)}(r) \leq \frac{1}{\sqrt{2}\pi^{(m+1)/2}} \Gamma\left(\frac{m}{2}\right) \overline{h}^{1/2} E(\|\xi_m\|) \left(\prod_{i=1}^{m} D_i^{1/2}\right) \sigma_{m-1} \int_0^r \frac{\rho^{m-1}}{(1+\rho^2)^{(m+1)/2}} \, d\rho, \tag{31}$$

where $\sigma_{m-1}$ is the $(m-1)$-dimensional area measure of $S^{m-1}$. The integral in the right-hand side is bounded by
$$\frac{\pi}{2} \left(\frac{r^2}{1+r^2}\right)^{(m-1)/2}.$$

Again using ($H_2$) and formula (9), we have the lower bound:

$$E(N^X) \geq \frac{\Gamma(m/2)}{\sqrt{2}\pi^{(m+1)/2}} \underline{h}^{1/2} E(\|\xi_m\|) \int_0^{+\infty} \left[\prod_{i=1}^m \left(\frac{D_i}{1+\rho^2} - \frac{E_i}{(1+\rho^2)^2}\right)^{1/2}\right] \frac{\rho^{m-1}}{(1+\rho^2)^{1/2}} \, d\rho$$

$$= \frac{\Gamma(m/2)}{\sqrt{2}\pi^{(m+1)/2}} \underline{h}^{1/2} E(\|\xi_m\|) \left(\prod_{i=1}^m D_i^{1/2}\right) \sigma_{m-1}$$

$$\times \int_0^{+\infty} \frac{\rho^{m-1}}{(1+\rho^2)^{(m+1)/2}} \prod_{i=1}^m \left(1 - \frac{F_i}{1+\rho^2}\right)^{1/2} d\rho$$

$$\geq \frac{\Gamma(m/2)}{\sqrt{2}\pi^{(m+1)/2}} \underline{h}^{1/2} E(\|\xi_m\|) \left(\prod_{i=1}^m D_i^{1/2}\right) \sigma_{m-1} \left(\frac{r_0^2 + 1/2}{r_0^2 + 1}\right)^{m/2}$$

$$\times \int_{\tau r_0}^{+\infty} \frac{\rho^{m-1}}{(1+\rho^2)^{(m+1)/2}} \, d\rho,$$

where $\tau$ has been chosen to satisfy (15).

To get a lower bound for the last integral, a direct integration by parts shows that:

$$\int_0^{+\infty} \frac{\rho^{m-1}}{(1+\rho^2)^{(m+1)/2}} \, d\rho > \frac{e^{-2}}{\sqrt{m}},$$

which implies

$$\int_{\tau r_0}^{+\infty} \frac{\rho^{m-1}}{(1+\rho^2)^{(m+1)/2}} \, d\rho > \frac{e^{-2}}{\sqrt{m}} - \frac{\pi}{2}\left(\frac{\tau^2 r_0^2}{1+\tau^2 r_0^2}\right)^{(m-1)/2}.$$

Then, choosing $m_0$ as required by condition (16), we get for $m \geq m_0$:

$$H_m^{(1)}(r_0) \leq C_1 \sqrt{m} \left(\frac{r_0^2}{r_0^2 + 1/2}\right)^{m/2} E(N^X), \tag{32}$$

where

$$C_1 = \pi e^2 \left(\frac{\bar{h}}{\underline{h}}\right)^{1/2} \frac{\sqrt{1+r_0^2}}{r_0}.$$

For the remainder, we must put together (27), (28), (30) and (32). One easily checks now that with the constants given by (17), the inequality of the statement holds true. □

### 3.1. Auxiliary lemma

**Lemma 1.** *Let $\gamma : \mathbb{R}^k \to \mathbb{R}$, $k \geq 1$ be defined as*

$$\gamma(c) = E(\|\xi + c\|),$$

where $\xi$ is a standard normal random vector in $\mathbb{R}^k$, and $c \in \mathbb{R}^k$ ($\|\cdot\|$ is the Euclidean norm in $\mathbb{R}^k$). Then

(i) $\gamma(0) = \sqrt{2}\frac{\Gamma((k+1)/2)}{\Gamma(k/2)}$.

(ii) $\gamma$ is a function of $\|c\|$ and verifies:

$$\gamma(c) \leq \gamma(0)\left(1 + \frac{1}{2k}\|c\|^2\right). \tag{33}$$

**Proof.** (i) In polar coordinates $\gamma(0) = \frac{\sigma_{j-1}}{(2\pi)^{j/2}} \int_0^{+\infty} \rho^j e^{-\rho^2/2}$. Using the change of variable $u = \rho^2$, we obtain the result.

(ii) That $\gamma$ is a function of $\|c\|$ is a consequence of the invariance of the distribution of $\xi$ under the isometries of $\mathbb{R}^k$. For $k = 1$, (33) follows from the exact computation

$$\gamma(c) = \sqrt{2/\pi}e^{-c^2/2} + c\int_{-c}^{c} \frac{1}{\sqrt{2\pi}}e^{-x^2/2}$$

and a Taylor expansion at $c = 0$, which gives

$$\gamma(c) \leq \sqrt{2/\pi}\left(1 + \frac{1}{2}c^2\right).$$

For $k \geq 2$, we write

$$\gamma(c) = E([(\xi_1 + a)^2 + \xi_2^2 + \cdots + \xi_k^2]^{1/2}) = G(a),$$

where $a = \|c\|$ and $\xi_1, \ldots, \xi_k$ are independent standard normal variables. Differentiating under the expectation sign, we get:

$$G'(a) = E\left(\frac{\xi_1 + a}{[(\xi_1 + a)^2 + \xi_2^2 + \cdots + \xi_k^2]^{1/2}}\right)$$

so that $G'(0) = 0$ due to the symmetry of the distribution of $\xi$.

One can differentiate formally once more, obtaining:

$$G''(a) = E\left(\frac{\xi_2^2 + \cdots + \xi_k^2}{[(\xi_1 + a)^2 + \xi_2^2 + \cdots + \xi_k^2]^{3/2}}\right). \tag{34}$$

For the validity of equality (34) for $k \geq 3$ one can use that if $d \geq 2$, $\frac{1}{\|x\|}$ is integrable in $\mathbb{R}^d$ with respect to the Gaussian standard measure. For $k = 2$ one must be more careful but it holds true. The other ingredient of the proof is that one can verify that $G''$ has a maximum at $a = 0$. Hence, on applying Taylor's formula, we get

$$G(a) \leq G(0) + \frac{1}{2}a^2 G''(0).$$

Check that $G''(0) = \frac{\sqrt{2}}{k} \frac{\Gamma((k+1)/2)}{\Gamma(k/2)}$ which, together with (i), gives:

$$\frac{G''(0)}{G(0)} = \frac{1}{k},$$

which implies (ii). □

## Acknowledgements

This work was supported by ECOS program U03E01. The authors thank Professor Felipe Cucker for useful discussions. Also, we thank two anonymous referees whose remarks have contributed to improving the paper.